\documentclass[11pt,final,journal,onecolumn,letterpaper]{IEEEtran}
\usepackage{amsmath,bm,bbm,amssymb}
\usepackage{undertilde}
\usepackage{cite}

\begin{document}
\title{Gibbs Random Fields and Markov Random Fields with Constraints}
\author{Levent~Onural
\thanks{L. Onural is with the Dept. of Electrical and Electronics Eng. of Bilkent University.}}

\bibliographystyle{IEEEtran}

\maketitle

\begin{abstract}
It was shown many times in the literature that a Markov random field is equivalent to a
Gibbs random field
when all realizations of the field have non-zero probabilities; the proofs are rather complicated.
A simpler proof, which is based directly on simple probability theory, is presented. Furthermore, it is
shown that the equivalence is still valid when there are constraints (zero probability realizations) of any type.
The equivalence extends to infinite size random fields, as well.
\end{abstract}



\section{Introduction}
\label{intro}
It is well known that a Gibbs random field has an equivalent Markov random field, and vice versa, 
provided that all realizations have nonzero probabilities; this requirement is known as the ``positivity condition''
\cite{hammersley71,besag1,spitzer,isham81}.
There are, however, many open questions when there are constraints, which in turn, impose zero probabilities to 
some outcomes (i.e., when there are impossible outcomes) 
\cite{derin89}.
Local constraints in the form of zero conditional probabilities are already incorporated into the theory
\cite{moussouris74}.
Some cases with global constraints are also demonstrated
\cite{onural87,onural91}.
The positivity condition is sufficient but not necessary; indeed, there are hints that it might not
be needed for the equivalence.
Here in this paper, the basic steps and the associated proofs to construct Gibbs or
Markov random fields are revisited.
Simpler alternative proofs are provided for many cases. 
It is also proven that Gibbs random fields and Markov random fields are equivalent even if there are impossible 
(zero probability, forbidden) outcomes.

\section{Preliminaries}
\label{prelim}

Let us start by the set of all outcomes ${\Omega}$, with its elements $a_i$ having all nonzero probabilities:
$P(a_i)\>=\>p_i\>>\>0$, for all $a_i$. 
Let us assume that these probabilities of elements, $P(a_i)\>=\>p_i$ may not be known, but all ratios,
$P(a_i)\>/P(a_j)\>=\>p_i/p_j\>=\>r_{ij}$ of such probabilities are known, for all $i$ and $j$.
It is easy to show that these ratios uniquely specify the unknown $P(a_i)$'s for all $i$.
To show this, let us choose an arbitrary $i$, and consider all ratios $p_j/p_i$, for all $j$.
Clearly,
\begin{equation}
\label{rat1}
\sum\limits_{j}\>\frac{p_j}{p_i}\>=\>\frac{1}{p_i}\quad,
\end{equation}
which implies
\begin{equation}
\label{rat2}
p_i\>=\>\frac{1}{\sum\limits_{j}\>\frac{p_j}{p_i}}\quad,
\end{equation}
and therefore,
\begin{equation}
\label{rat3}
p_j\>=\>\frac{p_i}{r_{ij}}\>=\>
\frac{1}{{r_{ij}\>}{\sum\limits_{j}\>\frac{p_j}{p_i}}}\quad\forall j\quad ,\quad \text {and for any $i$}\quad ,
\end{equation}
where the sum over $j$ indicates that the summation is running over all elements $a_j \in \Omega$.
Therefore, all unknown probabilities can be found from the known ratios of those probabilities.
The reason to explicitly include Eq.(\ref{rat3}) is to stress the fact that only one arbitrary element $a_i$
may be chosen as the reference, and all ratios of probabilities are given with respect to $p_i$ as $p_j/p_i$, 
for all $j$.
Indeed, it is sufficient to know only such a set of ratios of probabilities, instead of all ratios, to find
probabilities, $p_j$, for all $j$. 

Any constraint partitions the sample space $\Omega$ into two sets, ${\cal C}_k$, and its complement 
${\cal C}^c_k$.
There may be many such constraints, labeled by the index $k$, and in such a case, we are interested in the intersection ${\cal C}\>=\>\bigcap\limits_k\>{\cal C}_k$.
Therefore, any $a_i\in {\cal C}$ satisfies all constraints, and again, the sample space is partitioned into
${\cal C}$ and ${\cal C}^c$.

Now, let us concentrate on the conditional probabilities $P(a_i\>\big |\> {\cal C})$.
Simply, from the definition of the conditional probabilities, 
\begin{equation}
\label{cond}
P(a_i\>\big | \> {\cal C}) \>=\> 
\begin{cases}
\frac{P(a_i\cup {\cal C})}{P({\cal C})}\>=\>\frac{P(a_i)}{P({\cal C})} &\text{if $a_i\in {\cal C}$}\\
0&\text{else}
\end{cases}\quad .
\end{equation}
Clearly, $P({\cal C})\>>\>0$ provided that ${\cal C}$ is not the empty set, since all $a_i$'s have positive
probabilities.

Now, let us focus only on those outcomes in ${\cal C}$; by the way, conditional probabilities are probabilities, i.e.,
they satisfy all axioms and properties of a probability structure.
Therefore, we simply state that, i)
$P(a_i \>\big |\> {\cal C})\>>\>0$ for $a_i\in {\cal C}$, and, ii)
\begin{equation}
\label{condrat}
\frac{P(a_j\>\big |\> {\cal C})}{P(a_i\>\big |\> {\cal C})}\>=
\>\frac{\frac{P(a_j)}{P({\cal C})}}{\frac{P(a_i)}{P({\cal C})}}\>=\>\frac{P(a_j)}{P(a_i)}\>=\>\frac{p_j}{p_i}\>=\>r_{ij}
\quad .
\end{equation}
Therefore, given two outcomes where both satisfy the constraints, the ratio of their conditional 
probabilities is still
the same as the ratio of their (unconditional) probabilities, as expected.
This observation leads us to the conclusion that for any subset ${\cal C}$, if all the ratios of probabilities 
(or equivalently, conditional probabilities, given that both realizations are in ${\cal C}$), $r_{ij}$'s, are known, 
then the
conditional probabilities $P(a_j\>|\>{\cal C})$ for all $j\>\backepsilon a_j\in {\cal C}$ are also known as they
are induced by these ratios, following the same steps are in Eqs.(\ref{rat1}-\ref{rat3}); 
and furthermore, this will also induce
the (unconditional) probabilities of those elements, as well, if needed.

\section{GRF-MRF Equivalence in the Presence of Constraints}
\label{equiv}

Indeed, the equivalence of a GRF to its corresponding  MRF, and vice versa, is proven using these 
ratios of probabilities; the invariance of these ratios, whether for the unconditional case or the
conditional case, assures the GRF-MRF equivalence even if there are also constraints.  
And this is true for any constraint, as long as
the set of outcomes satisfying the constraints form a non-empty set.
As already stated,
constraints mean zero probabilities as also indicated by Eq.(\ref{cond}).  

Starting from the general, and simple, observations above, we can now turn our attention to GRFs and MRFs.
(We will use undertilde to describe random variables; no undertilde will be used for numerical values that
these random variables take.  We will use bold fonts to represent vectors (arrays).
For example, $P(\utilde{\bf x}\>=\> {\bf x})$ means the ``probability that the vector
random variable $\utilde{\bf x}$ takes the specific vector value ${\bf x}$.
We will also use the notation $P_{\utilde{\bf x}}({\bf x})$ for the same purpose, whenever
we feel this notation is more appropriate.
Indeed, whenever there is no ambiguity in the meaning, we will also use the shortened
notation $P({\bf x})$ to represent the same probability as described above. There is no specific meaning
attached to lower case or upper case symbols.)
As usual, we assume a set of indexed random variables, $\utilde{x}_l$; the number of elements in the 
set may be finite of infinite.
The index could be called the ``site'', but actually it may or may not be associated with a physical location.
The collection of all of those random variables $\utilde{x}_l$ for all $l$,
is a vector random variable $\utilde{\bf X}$; we will call  $\utilde{\bf X}$ as the ``random pattern''.
A realization of  $\utilde{\bf X}$, is denoted by  ${\bf X}$ which is a pattern over all sites.

As proven many times, conditional probabilities given for the MRF induce joint probabilities for all outcomes
$\utilde{\bf X}\>={\bf X}$, under the positivity condition 
\cite{hammersley71,besag1,spitzer,isham81,derin89}.
However, the provided proofs are unnecessarily complicated and lengthy.
Instead, a simple proof is a direct consequence of the trivial
discussion on the ratios of probabilites, as presented above:

Let us first prove that an MRF is also a GRF.
As before, the set ${\cal C}$ contains all outcomes that satisfy the constraints; therefore, the elements in
${\cal C}$ all have nonzero probabilities.
So, we can safely form the ratios of such probabilities, as,
\begin{equation}
\label{rat4}
\frac{P\left ( \utilde{\bf X}\>=\>{\bf X}_j \> \big |\> {\cal C} \right )}
     {P\left ( \utilde{\bf X}\>=\>{\bf X}_i \> \big |\> {\cal C} \right )}
\>=\> 
\frac{P\left ( \utilde{\bf X}\>=\>{\bf X}_j \right )}
     {P\left ( \utilde{\bf X}\>=\>{\bf X}_i \right )}
\>=\>
r_{ij}\quad\forall\quad i,j\quad .
\end{equation}
Let us choose patterns ${\bf X}$ from ${\cal C}$ 
such that their values at all sites, except the specific but arbitrary 
site $l$, are the same.
The number of such patterns is at least one; we are interested in cases where this number is greater than one since
otherwise the rest of the discussion is trivial. 
By the way, the number of such patterns could be quite small as a consequence of the constraints, and each such
distinct pattern has a different realization for $\utilde{x}_l$.
Assuming that there are more than one such patterns, we can write,
\begin{equation}
\label{mrf1}
P\left ( \utilde{\bf X}={\bf X}_j \> \big |\> {\cal C} \right )\>=\>
P\left ( \utilde{x}_l=v_j \> \big |\> \utilde{x}_m=x_m \,\forall m \,\backepsilon \, m\ne l, \,{\cal C}\right )
P\left (\utilde{x}_m=x_m \,\forall m \,\backepsilon \, m\ne l, \,{\cal C}\right ) \quad ,
\end{equation}
where $v_j$ is the realization (value) of the random variable $\utilde{x}_l$ at site $l$ within the pattern ${\bf X}_j$.
Therefore, the ratio in Eq.(\ref{rat4}) becomes,

\begin{equation}
\begin{aligned}
\label{mrf2}
\frac{P\left ( \utilde{\bf X}\>=\>{\bf X}_j \> \big |\> {\cal C} \right )}
     {P\left ( \utilde{\bf X}\>=\>{\bf X}_i \> \big |\> {\cal C} \right )}
\>&=\>
\frac
{P\left ( \utilde{x}_l\>=\>v_j \> \big |\> \utilde{x}_m\>=\>x_m \>\forall m \>\backepsilon \> m\ne l, \>{\cal C}\right )
P\left ( \utilde{x}_m\>=\>x_m \>\forall m \>\backepsilon \> m\ne l, \>{\cal C}\right )}
{P\left ( \utilde{x}_l\>=\>v_i \> \big |\> \utilde{x}_m\>=\>x_m \>\forall m \>\backepsilon \> m\ne l, \>{\cal C}\right )
P\left ( \utilde{x}_m\>=\>x_m \>\forall m \>\backepsilon \> m\ne l, \>{\cal C}\right )}
\nonumber\\
&=\>\frac
{P\left ( \utilde{x}_l\>=\>v_j \> \big |\> \utilde{x}_m\>=\>x_m \>\forall m \>\backepsilon \> m\ne l, \>{\cal C}\right )}
{P\left ( \utilde{x}_l\>=\>v_i \> \big |\> \utilde{x}_m\>=\>x_m \>\forall m \>\backepsilon \> m\ne l, \>{\cal C}\right )}\nonumber\\
\>&=\>\frac
{P\left ( \utilde{x}_l\>=\>v_j \> \big |\> \utilde{x}_m\>=\>x_m \> \text{for $m \in \eta_l$}, \>{\cal C}\right )}
{P\left ( \utilde{x}_l\>=\>v_i \> \big |\> \utilde{x}_m\>=\>x_m \> \text{for $m \in \eta_l$}, \>{\cal C}\right )}\nonumber\\
&\>=\>\frac
{P\left ( \utilde{x}_l\>=\>v_j \> \big |\> \utilde{x}_m\>=\>x_m \> \text{for $m \in \eta_l$}, \>{\cal C}\right )P({\cal C})}
{P\left ( \utilde{x}_l\>=\>v_i \> \big |\> \utilde{x}_m\>=\>x_m \> \text{for $m \in \eta_l$}, \>{\cal C}\right )P({\cal C})}\nonumber\\
\>&=\>\frac
{P\left ( \utilde{x}_l\>=\>v_j \> \big |\> \utilde{x}_m\>=\>x_m \> \text{for $m \in \eta_l$}\right )}
{P\left ( \utilde{x}_l\>=\>v_i \> \big |\> \utilde{x}_m\>=\>x_m \> \text{for $m \in \eta_l$}\right )}
\>=\>
r_{ij}\quad,
\end{aligned}
\end{equation}

\noindent
where we used the Markovianity, and $\eta_l$ is the neighborhood associated with site $l$ (see, for example,
\cite{derin89}, 
for the definition and the properties of the neighborhood).
Therefore, if all conditional probabilities,
\begin{equation}
\label{mrf3}
P\left ( \utilde{x}_l\>=\>v_l\> \big |\> \utilde{x}_m\>=\>x_m \> \text{for $m \in \eta_l$}, \>{\cal C}\right )\quad ,
\end{equation}
for all $l$ and for all allowed vaules $v_l$ for that 
location are known when there are constraints, or no constraints (i.e., when ${\cal C} = \Omega$, the sample space), 
and if the
Markovianity as indicated by 
\begin{equation}
\label{markov}
P\left ( \utilde{x}_l\>=\>v_l \> \big |\> \utilde{x}_m\>=\>x_m \>\forall m \>\backepsilon \> m\ne l\right )
\>=\>
P\left ( \utilde{x}_l\>=\>v_l \> \big |\> \utilde{x}_m\>=\>x_m \> \text{for $m \in \eta_l$} \right ) 
\end{equation}
holds for all $l$ and $v_l$ (i.e., if the field is a MRF), we can go backwards through the arguments and obtain
the result:
Given all conditional probabilities\linebreak 
$P\left ( \utilde{x}_l\>=\>v_l \> \big |\> \utilde{x}_m\>=\>x_m \> \text{for $m \in \eta_l$} \right )$ and the
set indicated by constraints ${\cal C}$, all ratios  $r_{ij}$ for all ${\bf X}_i$ and ${\bf X}_j$ in ${\cal C}$
are known; and from that (due to Eqs.(\ref{rat1}-\ref{rat3})) all conditional probabilities $P\left ( {\bf X}_i \> \big |\> 
{\cal C}\right )$ for all ${\bf X}_i \in {\cal C}$ are known and they are positive; 
and therefore, we can always write $P\left ( {\bf X}_i \>\big |\>{\cal C}\right )\>=\> k e^{-U({\bf X}_i)}$ due to
nonzero value of $P\left ( {\bf X}_i \>\big |\>{\cal C}\right )$, where $k$ is just a normalization constant to have
the sum of all probabilities equal to one, i.e.,
\begin{equation}
\label{normal1}
k\>=\> \sum\limits_{\text{all $i\>\backepsilon\> {\bf X}_i\in {\cal C}$}} e^{-U({\bf X}_i)} \quad .
\end{equation}
Furthermore, noting that,
\begin{equation}
\label{normal2}
P({\cal C})\>=\> \sum\limits_{\text{all $i\>\backepsilon\> {\bf X}_i\in {\cal C}$}} ke^{-U({\bf X}_i)} \quad ,
\end{equation}
we also know $P({\cal C})$, and therefore,
 $P({\bf X}_i)\>=\>P\left ({\bf X}_i \>\big |\> {\cal C}\right ) P({\cal C} )$ 
is also known, and is equal to $kP({\cal C} )  e^{-U({\bf X}_i)}$, for all ${\bf X}_i$, where $kP({\cal C} )$ is just the new
normalization constant to have the sum of all probabilities equal to one.
By the way, $U({\bf X})$ is called the ``energy'' of the pattern ${\bf X}$.
Therefore, a conditional (constraint satisfying) MRF is also a conditional GRF. 
The terms ``conditional Markov random fields (CMRF)'' and ``conditional Gibbs random fields (CGRF)''
are also used in
\cite{onural87,onural91}.
Please note that the proof presented above is also valid for the case where there are no constraints, i.e., when
${\cal C}\>=\Omega$, and therefore, we have also provided an alternative proof that a MRF is equivalent to a GRF
under the positivity condition, as well.
We believe that this proof is a lot simpler and straightforward than other known proofs, as given, for example 
in 
\cite{hammersley71,besag1,spitzer,isham81,derin89}.

Now we will prove that every conditional GRF is also a conditional MRF.
If $\utilde{\bf X}$ is a GRF, then each ${\bf X}_i\>\in\> {\cal C}$ has a probability in the form $ke^{-U({\bf X}_i)}$
as a consequence of the definition of a GRF;
we have shown above that the (unconditional) probabilities of these patterns are also in the form 
$kP({\cal C} )  e^{-U({\bf X}_i)}$, where $kP({\cal C} )$ is just the related normalization constant.
We can always decompose any function $U({\bf X}_i)$ into an additive form as,
\begin{equation}
\label{clpot}
U({\bf X}_i)\>=\>\sum\limits_{c\>\in\>Q}V_c({\bf X}_i)
\end{equation}
where $c$ is a ``clique'' which simply means a subset of indices (``sites'') $l$.
The set $Q$ is the set of all cliques.
The number of elements in a clique ranges from zero up to the maximum number of indices (sites) in a pattern; 
therefore, that maximum is infinity for infinite size patterns.
Therefore, the ratio of the probabilities of two patterns, ${\bf X}_j$ and ${\bf X}_i$ becomes,
\begin{equation}
\label{ratgibbs}
\frac{P\left ( \utilde{\bf X}\>=\>{\bf X}_j \>\big |\> {\cal C} \right )}
     {P\left ( \utilde{\bf X}\>=\>{\bf X}_i \>\big |\> {\cal C} \right )}
\>=\>
\frac{P\left ( \utilde{\bf X}\>=\>{\bf X}_j \right )}
     {P\left ( \utilde{\bf X}\>=\>{\bf X}_i \right )}
\>=\>
\frac{e^{-U({\bf X}_j)}}{e^{-U({\bf X}_i)}}
\>=\>
\frac{e^{-\sum\limits_{c\>\in\>Q}{V_c({\bf X}_j)}}}{e^{-\sum\limits_{c\>\in\>Q}V_c({\bf X}_i)}}\quad .
\end{equation}
Now, let us assume that the patterns ${\bf X}_i$, and ${\bf X}_j$ differ only at one specific site $l$, and they
have the same values at all other sites other than $l$.
In that case, the above ratio becomes,
\begin{equation}
\label{ratgibbs2}
\frac{e^{-\sum\limits_{c\>\in\>Q}{V_c({\bf X}_j)}}}{e^{-\sum\limits_{c\>\in\>Q}V_c({\bf X}_i)}}
\>=\>
\frac{\left ( e^{-\sum\limits_{c\>\in\>D_l}{V_c({\bf X}_j)}}\right ) \left ( e^{-\sum\limits_{c\>\in\>D_l^c}{V_c({\bf X}_j)}}\right )}
{\left ( e^{-\sum\limits_{c\>\in\>D_l}V_c({\bf X}_i)}\right )\left ( e^{-\sum\limits_{c\>\in\>D_l^c}V_c({\bf X}_i)}\right ) }\>
=\>
\frac{ e^{-\sum\limits_{c\>\in\>D_l}{V_c({\bf X}_j)}}}{ e^{-\sum\limits_{c\>\in\>D_l}V_c({\bf X}_i)}}\quad ,
\end{equation}
where $D_l$ is the set of only those cliques which have different realizations for ${\bf X}_i$, and ${\bf X}_j$ over
them, for a given $l$.  
In other words,
\begin{equation}
D_l\>=\>\left \{ c \> \>\big |\> {x}_{mi} \ne x_{mj} \> \text{for at least one $m \in c$}\right \}\quad ,
\end{equation}
where $x_{mi}$ and $x_{mj}$ are the realizations at site $m$ of the patterns ${\bf X}_i$, and ${\bf X}_j$,
respectively. 
Since  ${\bf X}_i$, and ${\bf X}_j$ may differ only at site $l$, as a consequence of the above assumption, 
an equivalent definition of $D_l$ can be given by
\begin{equation}
\label{set1}
D_l\>=\>\left \{ c\>\big | \> l\,\in\, c \right \}\quad .
\end{equation}
Please note that $D_l^c$ is the complement of $D_l$ in $Q$: $D_l^c\,=\,Q\setminus D_l$.
Noting, from Eq.(\ref{ratgibbs2}), that 
\begin{align}
\label{GRF1}
\frac{P\left ( \utilde{\bf X}={\bf X}_j \>\big |\> {\cal C} \right )}
     {P\left ( \utilde{\bf X}={\bf X}_i \>\big |\> {\cal C} \right )}
&=\>
\frac
{P\left ( \utilde{x}_l=v_j \>\big |\> \utilde{x}_m=x_m \>\forall m \>\backepsilon \> m\ne l, \>{\cal C}\right )}
{P\left ( \utilde{x}_l=v_i \>\big |\> \utilde{x}_m=x_m \>\forall m \>\backepsilon \> m\ne l, \>{\cal C}\right )}
\nonumber\\
&=\>
\frac
{\left ( e^{-\sum\limits_{c\>\in\>D_l}{V_c({\bf X}_j)}}\right )\left (e^{-\sum\limits_{c\>\in\>D_l^c}{V_c({\bf X}_j)}}\right )}
{\left ( e^{-\sum\limits_{c\>\in\>D_l}{V_c({\bf X}_i)}}\right )\left (e^{-\sum\limits_{c\>\in\>D_l^c}{V_c({\bf X}_i)}}\right )} 
\>=\>
\frac
{e^{-\sum\limits_{c\>\in\>D_l}{V_c({\bf X}_j)}}}{e^{-\sum\limits_{c\>\in\>D_l}V_c({\bf X}_i)}}\quad .
\end{align}
We conclude that $x_m$, $m\in\eta_l$ are sufficient to compute each term of the last expression above, to yield,
\begin{equation}
\label{markov3}
\frac
{P\left ( \utilde{x}_l\>=\>v_j \>\big |\> \utilde{x}_m\>=\>x_m \>\forall m \>\backepsilon \> m\ne l, \>{\cal C}\right )}
{P\left ( \utilde{x}_l\>=\>v_i \>\big |\>  \utilde{x}_m\>=\>x_m \>\forall m \>\backepsilon \> m\ne l, \>{\cal C}\right )}
\>=\>
\frac
{P\left ( \utilde{x}_l\>=\>v_j \>\big |\> \utilde{x}_m\>=\>x_m \> \text{for $m \in \eta_l$}, \>{\cal C}\right )}
{P\left ( \utilde{x}_l\>=\>v_i \>\big |\> \utilde{x}_m\>=\>x_m \> \text{for $m \in \eta_l$}, \>{\cal C}\right )}\quad ,
\end{equation}
and complete the proof that every GRF is also a MRF.
Please note that the proof is also valid if there are no constraints; i.e., when ${\cal C}\,=\,\Omega$.

\section{Observations,  Remarks and Conclusions}
\label{conc}

Based on the discussions above, we can conclude that,

\begin{itemize}
\item[*]
The only requirement to write a probability in exponential form is the
positivity of that probability; therefore, any field with positive probability realizations is a GRF.
\item[*]
A GRF has an equivalent MRF, and vice versa, when all the realizations (outcomes) have positive probabilities.
\item[*]
At this point, one could get the impression that MRFs, or equivalently GRFs, are mathematically so simple:
provided that outcomes (patterns) in a set of outcomes all have non-zero probabilities, the set forms a MRF (GRF).
This is true.
However, the MRF model becomes useful only when the neighborhoods $\eta_l$ are simple (small in size; i.e., have
few indices (sites) 
in it for every $l$); equivalently, GRF model becomes useful when the set of cliques $Q$ contains only
simple and few components.
An equivalent statement is that the clique potentials, $V_c({\bf X})$'s, are zero for most cliques (i.e., for
most subsets $c$ of sites).
In other words, the benefits prevail only when the direct statistical interactions among the random variables
at different sites are sparse.  
Obviously, every field with nonzero probability outcomes is a MRF when the sizes of the neighborhoods $\eta_l$'s
reach the size of the entire field; or equivalently, every field with nonzero probability outcomes is a GRF when
the sizes of $c$'s in $Q$ are allowed to reach the size of the entire field (i.e., $D^c\,=\,\varnothing$). 
In such cases, one may still call the field as a MRF with the neighborhood size equal to the size of the
entire field; or may choose to leave those extreme cases out of the definition and say that those cases are not 
MRFs (GRFs);
we prefer the first alternative in this paper. 
\item[*]
The neighborhood of a site is not necessarily near (here we assume that the sites form a lattice, and therefore,
a distance measure is applicable) to that site; indeed, the simplicity (sparsity) of a neighborhood scheme is
based on the number of sites in a neighborhood, and not where they are located.
\item[*]
Extension of the proofs for the Markov case where the ``interior'' is no longer a single 
site, but more than one site, is straightforward.
\item[*]
The discussion is valid both for causal fields, as well as non-causal ones; the causality is a direct
consequence of the neighborhood shape (see, for example 
\cite{won}
for the definition of ``causal'' and, ``non-causal''
random fields).
\item[*]
GRFs are usually simpler to handle than MRFs.
\item[*]
Inclusion of constraints, which in turn means zero probability outcomes, does not violate the MRF-GRF
equivalence or the basic structures of GRFs and the MRFs.
This is a consequence of the observation that the constraints result in a smaller set of non-zero probability 
realizations, and within that smaller set, all other features that leads to the equivalence is still valid and the
same as the unconstrained case.
\item[*]
These conclusions are valid both for the finite size or infinite size random fields.
\item[*]
Simple proofs are given for above statements; the proofs are direct consequences of elementary probability theory,
and techniques. 
\end{itemize}

The presented discussions clarify and provide answers to those long open questions related to constraints in MRFs or GRFs.

\end{document}